\documentclass[11pt]{article}
\usepackage{amssymb} 
\usepackage{amscd} 
\usepackage{amsthm} 
\usepackage{epsf}      %  replaces \usepackage{psfig} 
\newcommand{\stdspace}{\hskip 0.75em plus 0.15em\ignorespaces}
\newcommand{\co}{\colon\thinspace}    %  Colon with correct spacing for maps.
\newcommand{\nl}{\hfil\break}         %  New line.
\newcommand{\cl}{\centerline}         %  Centerline
\newcommand{\ppar}{\par\vskip 8pt plus4pt minus4pt}

\newcommand{\rk}[1]{\ppar{\bf #1}\stdspace}
\newtheoremstyle{plain}{14pt plus6pt minus6pt}{8pt plus3pt minus3pt}{\sl}%
{}{\bf}{}{0.75em}{\thmname{#1}\thmnumber{ #2}\thmnote{\rm\stdspace(#3)}}
\newtheoremstyle{definition}{14pt plus6pt minus6pt}{8pt plus3pt minus3pt}%
{\rm}{}{\bf}{}{0.75em}{\thmname{#1}\thmnumber{ #2}\thmnote{\sl\stdspace#3}}
\newtheoremstyle{remark}{14pt plus6pt minus6pt}{8pt plus3pt minus3pt}{\rm}%
{}{\bf}{}{0.75em}{\thmname{#1}\thmnumber{ #2}\thmnote{\sl\stdspace#3}}
\theoremstyle{plain}               
\def\sqr#1#2{{\vcenter{\vbox{\hrule  height.#2pt
	\hbox{\vrule width.#2pt height#1pt \kern#1pt \vrule width.#2pt}
	\hrule height.#2pt}}}}
\newtheorem{thm}{Theorem}
\newtheorem*{thmnn}{Theorem} 
\newtheorem{cor}{Corollary}[section]
\newtheorem{lem}[cor]{Lemma}

\theoremstyle{definition}
\newtheorem{defi}[cor]{Definition}

\newtheorem*{rem}{Remark}  
\makeatletter
\long\def\@caption#1[#2]#3{\par\addcontentsline{\csname
  ext@#1\endcsname}{#1}{\protect\numberline{\csname
  the#1\endcsname}{\ignorespaces #2}}\begingroup
    \@parboxrestore
    \small
    \@makecaption{\csname fnum@#1\endcsname}{\ignorespaces #3}\par
  \endgroup}
\newenvironment{prf}[1][\proofname]{\par
  \normalfont
  \topsep6\p@\@plus6\p@ \trivlist
  \item[\hskip\labelsep\bf
    #1]\ignorespaces
}{%
  \qed\endtrivlist{\parskip 12pt plus 6pt minus 6pt \par}
}
\makeatother
\newcommand{\Z}{{\mathbb{Z}}}

\newcommand{{\R}}{{\mathbb{R}}}

\newcommand{\hra}{\hookrightarrow}

\newcommand{\ra}{\longrightarrow}

\newcommand{\gropes}{
$$\begin{picture}(140,107) \small
    \put(-107,-10)       {\epsfbox{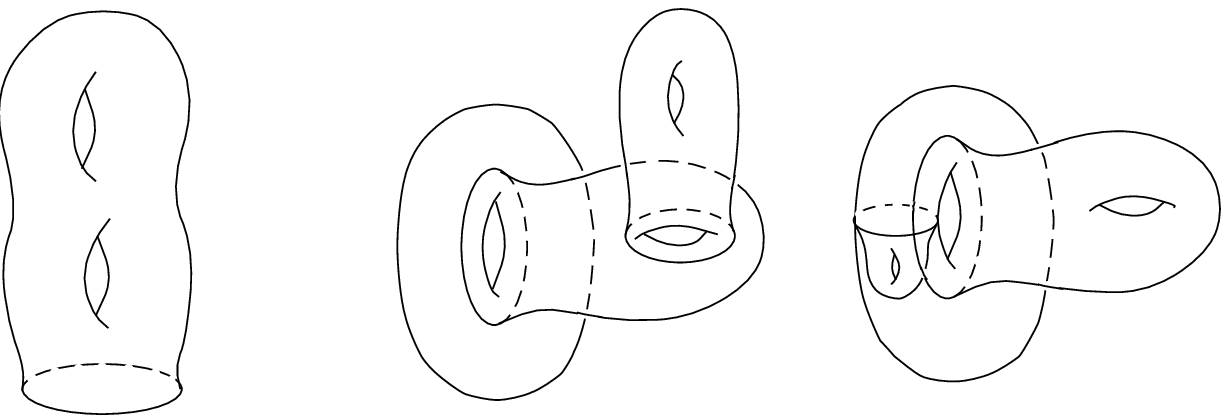}}
    \end{picture}$$
}
\newcommand{\gropestrees}{
$$\begin{picture}(140,150) \small
    \put(-85,-10)         {\epsfbox{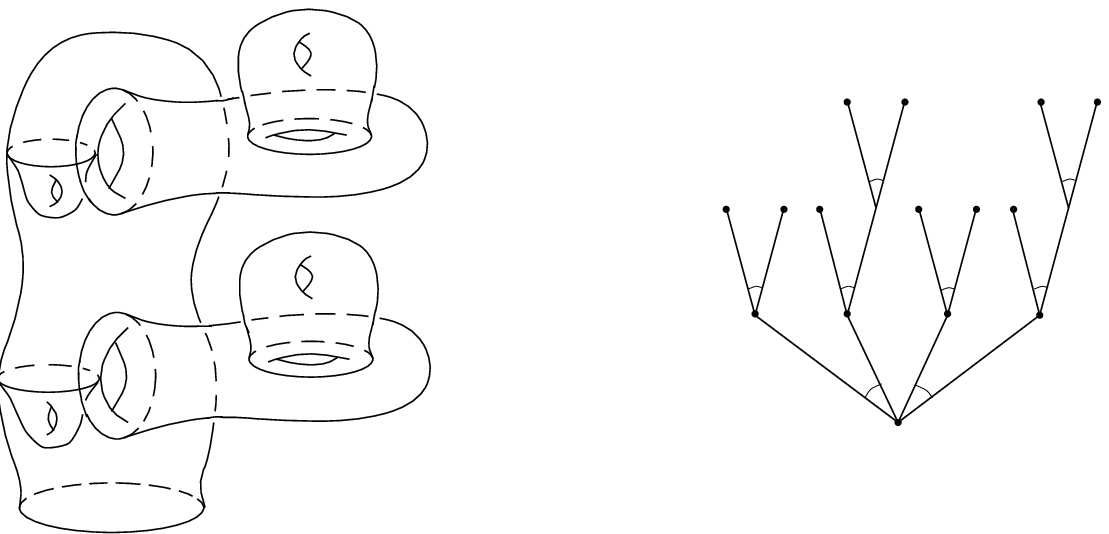}}
    \put(167,10)         {root}
    \put(162,125)        {leaves}
    \end{picture}$$
}
\newcommand{\torus}{
$$\begin{picture}(140,165) \small
    \put(-92,-10)       {\epsfbox{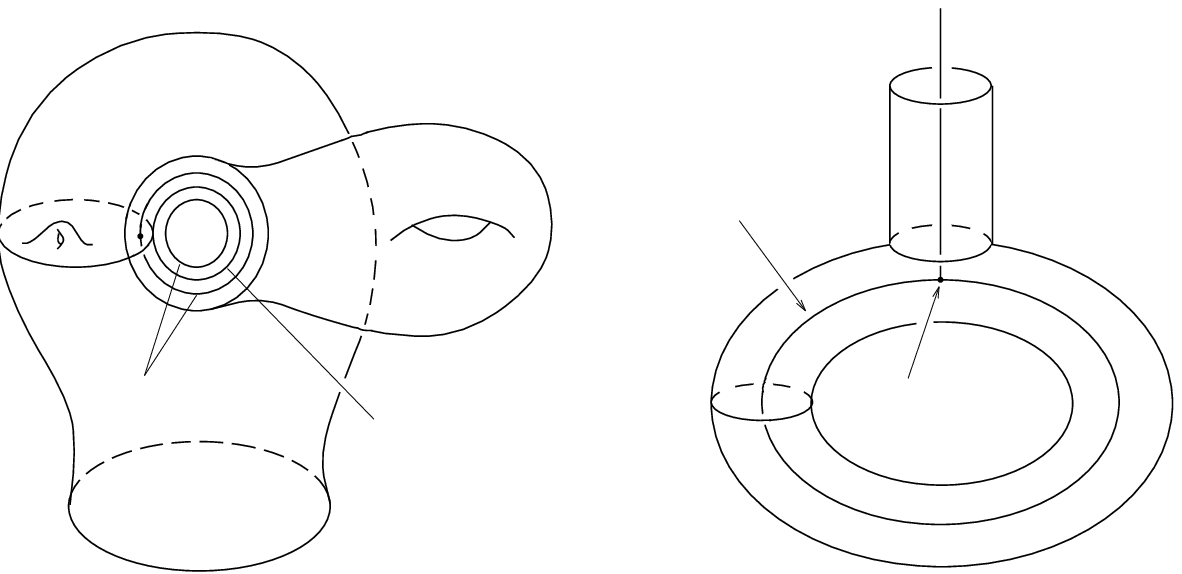}}
    \put(-93,35)       {${\Sigma}_i$}
    \put(-112,85)      {${\beta}_{i,j}$}
    \put(-53,114)      {${\gamma}={\alpha}_{i,j}$}
    \put(-63,39)       {$A_{i,j}$}
    \put(10,27)        {${\alpha}_{i,j}$(displaced)}
    \put(26,65)        {$X'$}
    \put(87,39)        {$m_{i-1}$}
    \put(90,95)        {${\beta}_{i-1,n}\subset{\Sigma}_{i-1}$}
    \put(183,145)      {${\Sigma}_i$}
    \put(199,128)      {$m_i$}
    \put(160,39)       {${\alpha}_{i-1,n}$}
    \end{picture}$$
}
\newcommand{\dualtree}{
$$\begin{picture}(140,120) \small
    \put(-106,-10)       {\epsfbox{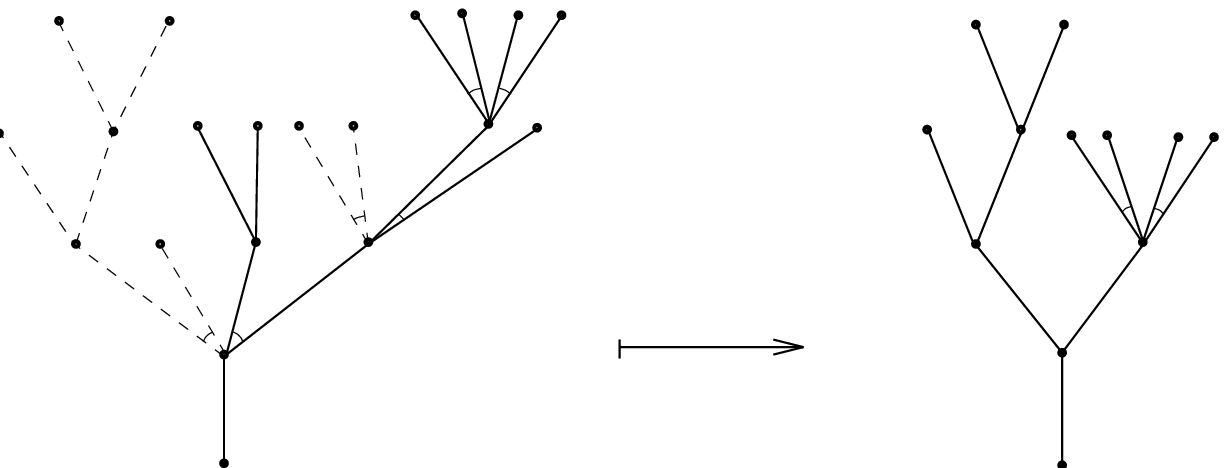}}
    \put(-32,-5)      {$T_{X}$}
    \put(211,-5)      {$T_{G}$}
    \put(2,50)        {${\alpha}_{i-1,n}$}
    \put(12,88)       {${\alpha}_{i,j}$}
    \put(51,87)       {${\beta}_{i,j}$}
    \put(146,87)      {$m_1$}
    \put(158,53)      {$m_2$}
\end{picture}$$
}
\newcommand{\epsilonfigure}{
$$\begin{picture}(160,125) \small
    \put(-85,0)        {\epsfbox{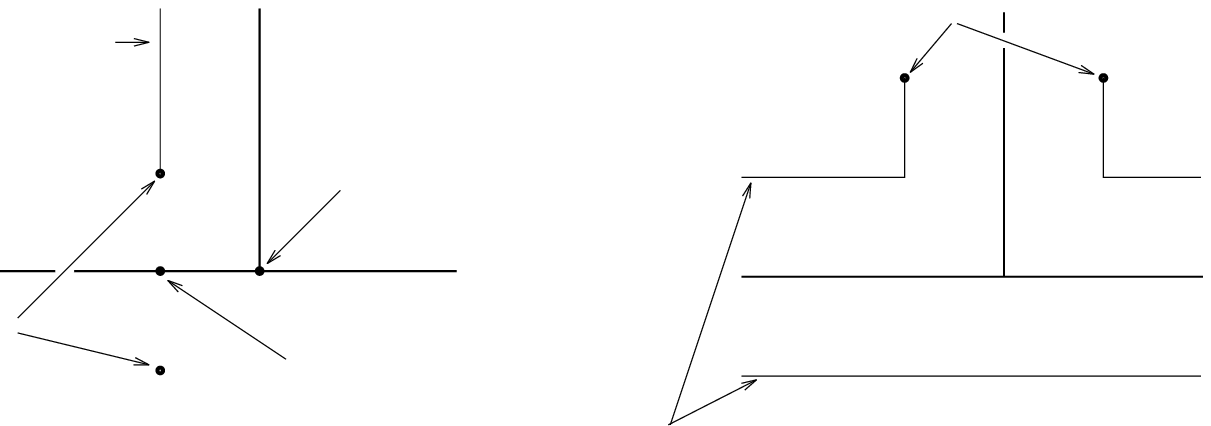}}
    \put(-100,28)     {$A_{i,j}$}
    \put(-103,110)    {a parallel}
    \put(-103,98)     {copy of $X'$}
    \put(-5,110)      {$X'$}
    \put(13,71)       {${\alpha}_{i,j}$}
    \put(48,42)       {${\beta}_{i,j}\subset {\Sigma}_i$}
    \put(75,-9)       {punctured $B_{i,j}$}
    \put(185,120)     {$m_i$}
    \put(-2,13)       {${\alpha}_{i,j}$(displaced)}
    \put(-28,85)      {$\epsilon$}
    \put(-42,57)      {$\epsilon$}
    \put(187,85)      {$\epsilon$}
    \put(160,56)      {$\epsilon$}
\end{picture}$$
}
\newcommand{\genusone}{
$$\begin{picture}(140,115) \small
    \put(-70,-10)       {\epsfbox{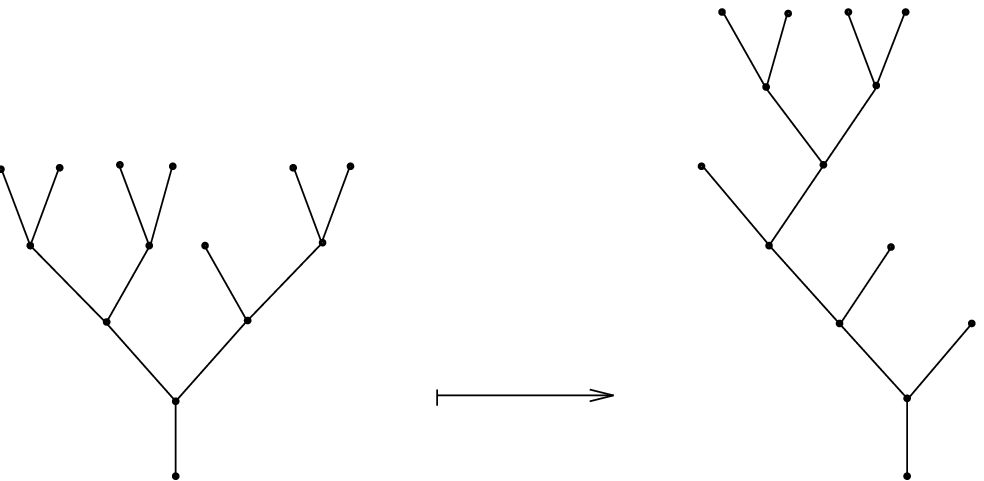}}
    \put(-10,-5)       {$T_{X}$}
    \put(200,-5)       {$T_{G}$}
    \put(3,79)         {$\gamma$}
    \put(113,79)       {$m_1$}
    \put(132,55)       {$m_2$}
    \put(152,32)       {$m_3$}
\end{picture}$$
}
\newcommand{\composition}{
$$\begin{picture}(100,155) \small
    \put(-75,-5)       {\epsfbox{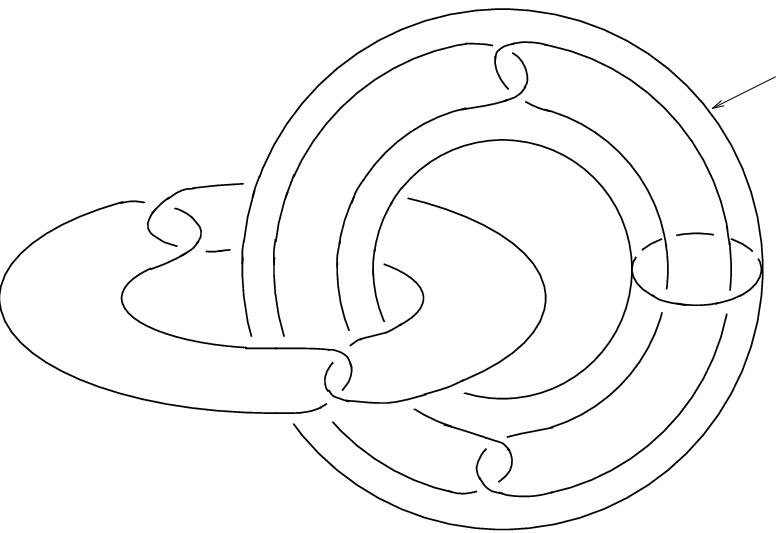}}
    \put(-52,97)      {$L$}
    \put(135,20)      {${\phi}(S^1\times D^2)$}
    \put(150,70)      {${\phi}(\wedge)$}
    \put(151,128)     {$\wedge'$}
    \end{picture}$$
}
\begin{document}
\setlength{\abovedisplayskip}{6pt plus3pt minus3pt}   %  Reduces the space
\setlength{\belowdisplayskip}{6pt plus3pt minus3pt}   %  around displays.
%
%     Definitions for title page:
%
\def\title#1{\def\thetitle{#1}}
\def\authors#1{\def\theauthors{#1}}
\def\address#1{\def\theaddress{#1}}
\def\email#1{\def\theemail{#1}}
\def\url#1{\def\theurl{#1}}
\long\def\abstract#1\endabstract{\long\def\theabstract{#1}}
\def\primaryclass#1{\def\theprimaryclass{#1}}
\def\secondaryclass#1{\def\thesecondaryclass{#1}}
\def\keywords#1{\def\thekeywords{#1}}
%
%                          Title page  
%
% This page will be reformatted by the Journal so you can use any
% format you like.  In particular you may use the standard latex  
% format for heading material if you wish.  
%
% Acknowledgements should not appear on this page. 
% Place these at the end of your introduction.
%  
\input gtoutput
\volumenumber{1}\papernumber{5}\volumeyear{1997}
\published{26 October 1997}
\pagenumbers{51}{69}

\title{Alexander Duality, Gropes and Link Homotopy}     
\shorttitle{Alexander duality, gropes and link homotopy}

\authors{Vyacheslav S Krushkal\\Peter Teichner}     

\address{Department of Mathematics, Michigan State University\\                    %address
East Lansing, MI 48824-1027, USA\\\vspace{.05cm}\\
{\rm Current address:}\ \ Max-Planck-Institut f\"{u}r Mathematik\\
Gottfried-Claren-Strasse 26, D-53225 Bonn, Germany
\\\vspace{.05cm}\\\cl{\rm and}\\\vspace{.05cm}\\
Department of Mathematics\\University of California in San Diego\\
La Jolla, CA, 92093-0112, USA}
\email{krushkal@math.msu.edu\\teichner@euclid.ucsd.edu}
\let\theurl\relax

\abstract 

We prove a geometric refinement of Alexander duality for certain 
$2$--complexes, the so-called {\em gropes}, embedded into $4$--space. 
This refinement can be roughly formulated as saying that
$4$--dimensional Alexander duality preserves the {\em disjoint
Dwyer filtration.} 

In addition, we give new proofs and extended versions of two lemmas  
of Freedman and Lin which are of central importance in the 
{\em A-B--slice problem}, the main open problem in the classification 
theory of topological $4$--manifolds. Our methods are group theoretical, 
rather than using Massey products and Milnor $\mu$--invariants as in 
the original proofs.
\endabstract

\primaryclass{55M05, 57M25}
\secondaryclass{57M05, 57N13, 57N70}


\proposed{Robion Kirby}\received{17 June 1997}
\seconded{Michael Freedman, Ronald Stern}\revised{17 October 1997}

\maketitlepage
%%%%%%%%%%%%%%%%%%%%   End of title page
%
%%%%%%%%%%%%%%%%%%%    Start of main body of article

\section{Introduction} \label{introduction}
Consider a finite complex $X$ PL--embedded into the $n$--dimensional sphere $S^n$. Alexander duality identifies the
(reduced integer) homology $H_i(S^n\smallsetminus X)$ with the cohomology $H^{n-1-i}(X)$. This
implies that the homology (or even the stable homotopy type) of the complement cannot distinguish between possibly different embeddings of $X$ into $S^n$.
Note that there cannot be a duality for homotopy groups as one can see by considering the fundamental group of
classical knot complements, ie the case $X=S^1$ and $n=3$.

However, one can still ask whether additional information about $X$ does lead to additional information about $S^n
\smallsetminus X$. For example, if $X$ is a smooth closed $(n-1-i)$--dimensional manifold then the
cohomological fundamental class is dual to a {\em spherical} class in $H_i(S^n\smallsetminus X)$. Namely, it is
represented by any {\em meridional} $i$--sphere which by definition is the boundary of a normal disk at a point in $X$.
This geometric picture explains the dimension shift in the Alexander duality theorem.

By reversing the roles of $X$ and $S^n \smallsetminus X$ in this example we see that it is {\em not} true that $H_i(X)$
being spherical implies that $H_{n-1-i}(S^n \smallsetminus X)$ is spherical. However, the following result shows that
there is some kind of improved duality if one does {\em not} consider linking dimensions. One should think of the
Gropes in our theorem as means of measuring how spherical a homology class is.

\begin{thm} \label{duality}\rm (Grope Duality)\ \ \sl
If $X \subset S^4$ is the disjoint union of closed embedded Gropes of class~$k$ then $H_2(S^4\smallsetminus X)$ is
freely generated by $r$ disjointly embedded closed Gropes of class~$k$. Here $r$ is the rank of $H_1(X)$.
Moreover, $H_2(S^4\smallsetminus X)$ cannot be generated by $r$ disjoint
maps of closed gropes of class $k+1$.
\end{thm}
As a corollary to this result we show in \ref{grope concordance} that certain Milnor $\mu$--invariants of a link in
$S^3$ are unchanged under a {\em Grope concordance}.

\begin{figure}[ht]
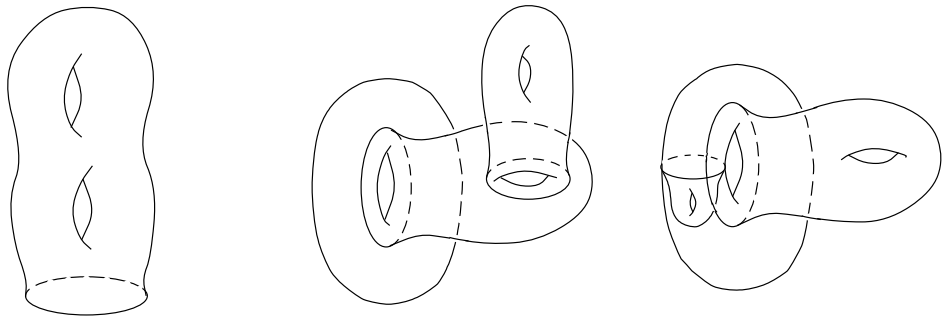
 
\gropes
\caption{A grope of class 2 is a surface -- two closed gropes of class 4}
\label{gropes}
\end{figure}

The {\em Gropes} above are framed thickenings of very simple 2--complexes, called {\em gropes}, which are inductively built out of surface stages, see Figure~\ref{gropes} and Section~\ref{facts}. 
For example, a grope of class~$2$ is just a surface with a single boundary component and gropes of bigger class contain information about the lower central series of the fundamental group. Moreover, every closed grope has a fundamental class in $H_2(X)$ and one obtains a geometric definition of the {\em Dwyer filtration}
\[
\pi_2(X) \subseteq \dots \subseteq \phi_k(X) \subseteq \dots \subseteq\phi_3(X)\subseteq\phi_2(X)=H_2(X)
\]
by defining $\phi_k(X)$ to be the set of all homology classes represented
by maps of closed gropes of class~$k$ into $X$.
Theorem \ref{duality} can thus be roughly formulated as saying that  
$4$--dimensional Alexander duality preserves the {\em disjoint
Dwyer filtration.} 

Figure~\ref{gropes} shows that each grope has a certain ``type'' which measures how the surface stages are attached. In Section~\ref{facts} this will be made precise using certain rooted trees, compare Figure~\ref{gropes and trees}. In Section~\ref{dualitysection}
we give a simple algorithm for obtaining the trees corresponding to the dual Gropes constructed in
Theorem~\ref{duality}.

The simplest application of Theorem~\ref{duality} (with class $k=2$) is as follows. Consider the standard embedding of
the $2$--torus $T^2$ into $S^4$ (which factors through the usual unknotted picture of $T^2$ in $S^3$). Then the boundary of
the normal bundle of $T^2$ restricted to the two essential circles gives two disjointly embedded tori representing
generators of $H_2(S^4 \smallsetminus T^2) \cong \Z^2$. Since both of these tori may be surgered to (embedded) spheres,
$H_2(S^4 \smallsetminus T^2)$ is in fact spherical. 
However, it cannot be generated by two maps of $2$--spheres with
{\em disjoint} images, since a map of a sphere may be replaced
by a map of a grope of arbitrarily big class.

This issue of disjointness leads us to study
the relation of gropes to classical {\em link homotopy}. We use Milnor group techniques to give new proofs and improved versions of the two central results of \cite{FL}, namely the {\em Grope Lemma} and the {\em Link Composition Lemma}.
Our generalization of the grope lemma reads as follows.

\begin{thm} \label{grope lemma} \sl
Two $n$--component links in $S^3$ are link homotopic if and only if they cobound disjointly immersed annulus-like gropes
of class~$n$ in $S^3 \times I$.
\end{thm}

This result is stronger than the version given in \cite{FL} where the authors only make a comparison with the trivial
link. Moreover, our new proof is considerably  shorter than the original one. 

Our generalization of the link composition lemma is formulated as Theorem~\ref{composition} in
Section~\ref{compositionsection}. The reader should be cautious about the proof given in \cite{FL}. It turns out that
our Milnor group approach contributes a beautiful feature to Milnor's algebraization of link homotopy: He proved in
\cite{Milnor1} that by forgetting one component of the unlink one gets an abelian normal subgroup of the Milnor group which is the additive group of a certain ring $R$. We observe that the {\em Magnus expansion} of the free Milnor groups arises naturally from considering the conjugation action of the quotient group on this ring $R$. Moreover, we
show in Lemma~\ref{multiplication} that ``composing'' one link into another corresponds to multiplication in that particular
ring $R$. This fact is the key in our proof of the link composition lemma.

Our proofs completely avoid the use of Massey products and Milnor $\! \mu$--invariants and we feel that they are more
geometric and elementary than the original proofs. This might be of some use in studying the still unsolved {\em
A-B--slice problem} which is the main motivation behind trying to relate gropes, their duality and link homotopy. It is
one form of the question whether topological surgery and s--cobordism theorems hold in dimension~$4$ without
fundamental group restrictions. See \cite{FT1} for new developments in that area.

\noindent
{\em Acknowledgements:}  It is a pleasure to thank Mike Freedman for many important discussions and for providing an
inspiring atmosphere in his seminars. In particular, we would like to point out that the main construction of
Theorem~\ref{duality} is reminiscent of the methods used in the {\em linear grope height raising} procedure of
\cite{FT2}.
The second author would like to thank the Miller foundation at 
UC Berkeley for their support.

\section{Preliminary facts about gropes and the lower \nl central series} \label{facts}
\begin{figure}[ht]
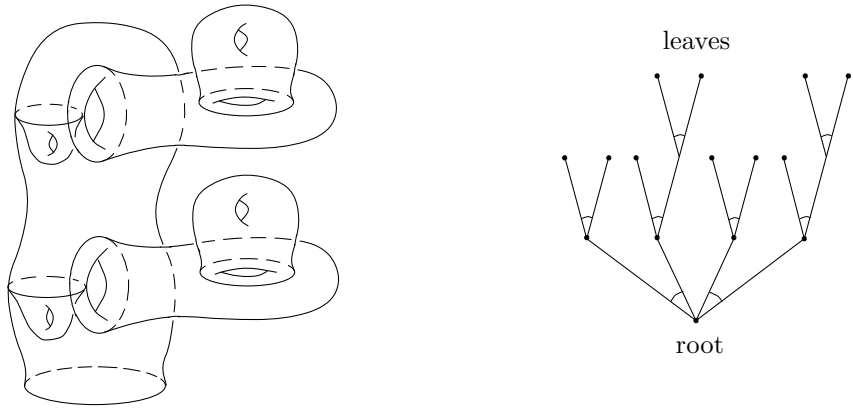
 
\gropestrees
\caption{A grope of class $5$ and the associated tree}
\label{gropes and trees}
\end{figure}

The following definitions are taken from \cite{FT2}.
\begin{defi} \rm
  A {\it grope} is a
special pair (2--complex, circle).  A grope has a {\em class}
$k=1, 2,\dots,
\infty$.  For $k=1$  a grope is defined to be the pair
(circle, circle).  For $k=2$ a grope is precisely a
compact oriented surface $\Sigma$ with a single boundary
component.  For $k$ finite a $k$--{\it grope} is
defined inductively as follow:  Let $\{\alpha_i,
\beta_i, i=1, \dots, {\rm genus}\}$ be a standard
symplectic basis of circles for $\Sigma$.  For any
positive integers $p_i, q_i$ with $p_i+q_i\ge k$ and
$p_{i_0} + q_{i_0} = k$ for at least one index $i_0$,
a
$k$--grope is formed by gluing $p_i$--gropes to each
$\alpha_i$ and
$q_i$--gropes to each $\beta_i$.
\end{defi}
The important information about the ``branching" of a
grope can be very well captured in a rooted tree as
follows:  For $k=1$ this tree consists of a single
vertex $v_0$ which is called the {\em root}.  For $k=2$ one adds
$2\cdot$genus$(\Sigma)$ edges to $v_0$ and may
label the new vertices by $\alpha_i, \beta_i$.
Inductively, one gets the tree for a $k$--grope which
is obtained by attaching $p_i$--gropes to $\alpha_i$
and $q_i$--gropes to $\beta_i$ by identifying the
roots of the $p_i$--(respectively $q_i$--)gropes with the
vertices labeled by $\alpha_i$(respectively $\beta_i$).
Figure~\ref{gropes and trees} below should explain the
correspondence between gropes and trees.

 Note that the vertices of the tree which are above the root
$v_0$ come in pairs corresponding to the symplectic
pairs of circles in a surface stage and that such
rooted paired trees correspond bijectively to
gropes.  Under this bijection, the {\it leaves}
($:=$ 1--valent vertices) of the tree correspond to
circles on the grope which freely generate its
fundamental group.  We will sometimes refer to these
circles as the {\it tips} of the grope.  The boundary of
the first stage surface $\Sigma$ will be referred to
as the {\it bottom} of the grope.

Given a group $\Gamma$, we will denote by $\Gamma^k$ the $k$-th
term in the lower central series of $\Gamma$, defined inductively by
$\Gamma^1:=\Gamma$ and $\Gamma^k:=[\Gamma,\Gamma^{k-1}]$, the characteristic subgroup of $k$--fold commutators in $\Gamma$.

\begin{lem}\label{maximal k} \sl
{\rm (Algebraic interpretation of gropes \cite[2.1]{FT2})}\ \ 
For a space X, a loop $\gamma$ lies in
$\pi_1(X)^k, 1 \le k < \omega$, if and only if $\gamma$
bounds a map of some $k$--grope.  Moreover,
the class of a grope $(G,\gamma)$ is the maximal $k$
such that $\gamma \in \pi_1(G)^k$.
\end{lem}

A {\it closed} $k$--grope is a 2--complex made by replacing a 2--cell in
$S^2$ with a $k$--grope.
A closed grope is sometimes also called a {\em sphere-like} grope. Similarly, one has
{\em annulus-like} $k$--gropes which are obtained from an annulus by replacing a 2--cell with a $k$--grope.
Given a space $X$, the Dwyer's subgroup $\phi_k(X)$
of $H_2(X)$ is the set of all homology classes represented
by maps of closed gropes of class~$k$ into $X$. Compare \cite[2.3]{FT2} for a translation to Dwyer's original
definition.

\begin{thmnn} \rm(Dwyer's Theorem \cite{Dwyer})\ \ \sl
Let $k$ be a positive integer and let
$f\co$ $X\longrightarrow Y$ be a map inducing an isomorphism
on $H_1$ and an epimorphism on
$H_2/\phi_k$. Then $f$ induces an
isomorphism on $\pi_1/(\pi_1)^k$.
\end{thmnn}

A {\em Grope} is a special
``untwisted'' $4$--dimensional thickening of a grope $(G,\gamma)$;
it has a preferred solid torus (around the base
circle $\gamma$) in its boundary. This ``untwisted'' thickening
is obtained by first embedding $G$ in ${\R}^3$ and taking its
thickening there, and then crossing it with the interval $[0,1]$.
The definition of a Grope is independent of the
chosen embedding of $G$ in ${\R}^3$. 
One can alternatively define it by a thickening of $G$ such that all relevant relative Euler numbers vanish.
Similarly, one defines sphere- and annulus-like Gropes, the capital letter
indicating that one should take a $4$--dimensional untwisted thickening of the corresponding 2--complex.

\section{The Grope Lemma} \label{gropesection}
We first recall some material from \cite{Milnor1}.
Two $n$--component links $L$ and $L'$
in $S^3$ are said to be {\em link-homotopic}
if they are connected by a 1--parameter family of immersions
such that distinct components stay disjoint at all times.
$L$ is said to be {\em homotopically trivial} if it is
link-homotopic to the unlink. $L$ is {\em almost homotopically
trivial} if each proper sublink of $L$ is homotopically trivial.

For a group $\pi$ normally generated by $g_1,\ldots,g_k$ its
{\em Milnor group} $M\pi$ (with respect to $g_1,\ldots,g_k$)  is
defined to be the quotient of $\pi$ by the normal subgroup
generated by the elements $[g_i,g_i^h]$, where $h\in\pi$ is arbitrary.
Here we use the conventions
\[ [g_1,g_2]:=g_1\cdot g_2 \cdot g_1^{-1}\cdot g_2^{-1} \; 
{\rm and } \; g^h:=h^{-1} \cdot g \cdot h.
\]
$M\pi$ is nilpotent of class~$\leq k+1$, ie it is
a quotient of $\pi/(\pi)^{k+1}$, and is generated by the quotient images
of $g_1,\ldots, g_k$, see \cite{FT1}.
The Milnor group $M(L)$ of a link $L$ is defined to be
$M\pi_1(S^3\smallsetminus L)$ with respect to its meridians $m_i$.
It is the largest common quotient of the fundamental
groups of all links link-homotopic to $L$, hence one obtains:

\begin{thmnn}\rm (Invariance under link homotopy \cite{Milnor1})\ \ \sl
If $L$ and $L'$ are link homotopic
then their Milnor groups are isomorphic.
\end{thmnn}
The track of a link homotopy in $S^3 \times I$ gives disjointly immersed annuli with the additional property of being
mapped in a level preserving way. However, this is not really necessary for $L$ and $L'$ to be link homotopic, as the
following result shows.

\begin{lem}\rm (Singular concordance implies homotopy \cite{Giffen}, 
\cite{Goldsmith}, \cite{Lin})\ \ \label{concordance} \sl 
\nl
If $L\subset S^3\times\{ 0\}$ and $L'\subset S^3\times\{ 1\}$
are connected in $S^3\times I$ by disjointly immersed annuli then $L$ and $L'$
are link-homotopic.
\end{lem}
\begin{rem} \rm
This result was recently generalized to all dimensions, see \cite{Teichner}.
\end{rem}
Our Grope Lemma (Theorem~\ref{grope lemma} in the introduction) further weakens the conditions on the objects that
connect $L$ and $L'$.

\rk{Proof of Theorem~\ref{grope lemma}}
Let $G_1,\ldots,G_n$ be disjointly immersed annulus-like gropes of class $n$ connecting $L$ and $L'$ in $S^3
\times I$. To apply the above Lemma~\ref{concordance}, we want to replace one $G_i$ at a time by an immersed annulus
$A_i$ in the complement of all gropes and annuli previously constructed.

Let's start with $G_1$. Consider the circle $c_1$ which consists of the union of the first component $l_1$ of $L$, then
an arc in $G_1$ leading from $l_1$ to $l_1'$, then the first component $l_1'$ of $L'$ and finally the same arc back to
the base point. Then the $n$--grope $G_1$ bounds $c_1$ and thus $c_1$ lies in the $n$-th term of the lower central
series of the group $\pi_1(S^3 \times I \smallsetminus G)$, where $G$ denotes the union of $G_2,\ldots,G_n$.  As first
observed by Casson, one may do finitely many finger moves on the bottom stage surfaces of G (keeping the components $G_i$ disjoint) such that the
natural projection induces an isomorphism
\[ \pi_1(S^3 \times I \smallsetminus G) \cong M\pi_1(S^3 \times I \smallsetminus G) \]
(see \cite{FT1} for the precise argument, the key idea being that the relation $[m_i,m_i^h]$ can be achieved by a self
finger move on $G_i$ which follows the loop $h$.) But the latter Milnor group is normally generated by $(n-1)$
meridians and is thus nilpotent of class~$\leq n$. In particular, $c_i$ bounds a disk in $S^3 \times I \smallsetminus
G$ which is equivalent to saying that $l_1$ and $l_1'$ cobound an annulus $A_1$, disjoint from $G_2,\ldots,G_n$.

Since finger moves only change the immersions and not the type of a 2--complex, ie an immersed annulus stays an
immersed annulus, the above argument can be repeated $n$~times to get disjointly immersed annuli $A_1,\ldots,A_n$ connecting $L$
and $L'$.
\qed

\section{Grope Duality} \label{dualitysection}
\begin{figure}[ht]
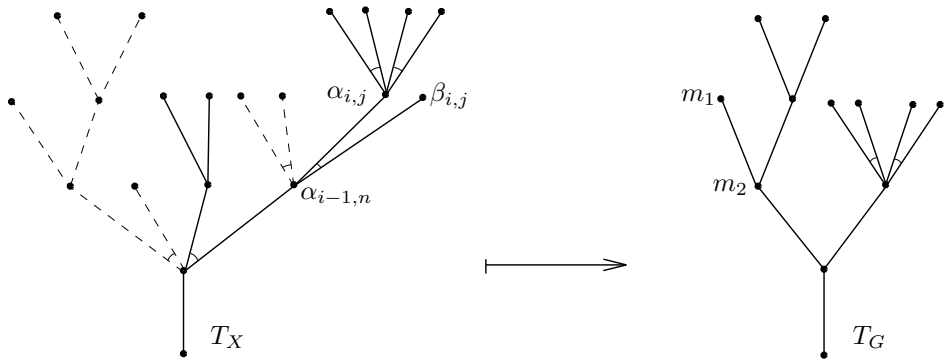
 
\dualtree
\caption{A dual tree: The branches in $T_{X}$ to be erased
are drawn with dashed lines.} \label{dualtree}
\end{figure}

In this section we give the proof of Theorem~\ref{duality} and a 
refinement which explains what the trees corresponding to the dual 
Gropes look like. Since we now consider closed gropes, the following variation of the correspondence to trees turns out to be extremely useful. 
Let $G$ be a closed grope and let $G'$ denote $G$ with a small $2$--cell removed from its bottom stage. 
We define the tree $T_G$ to be the tree corresponding to $G'$ (as defined in Section~\ref{facts}) together with an edge added to the root vertex. This edge represents the deleted $2$--cell and it turns out to be useful to define the root of $T_G$ to be the $1$--valent vertex of this new edge. See Figure~\ref{dualtree} for an example of such a tree.

\rk{Proof of Theorem~\ref{duality}}
Abusing notation, we denote by $X$ the core grope of the given $4$--dimensional Grope in $S^4$. Thus $X$ is a $2$--complex which has a particularly simple thickening in $S^4$ which we may use as a regular neighborhood.
All constructions will take place in this regular neighborhood, so we may assume that $X$ has just one connected component.
Let $\{\alpha_{i,j},\beta_{i,j}\}$ denote a standard symplectic basis
of curves for the $i$-th stage $X_i$ of $X$; these curves correspond
to vertices at a distance $i+1$ from the root in the associated tree. Here $X_1$ is the bottom stage and thus a closed connected surface. For $i>1$, the $X_i$ are disjoint unions of punctured surfaces.
They are attached along some of the curves $\alpha_{i-1,j}$ or $\beta_{i-1,j}$.

\begin{figure}[ht]
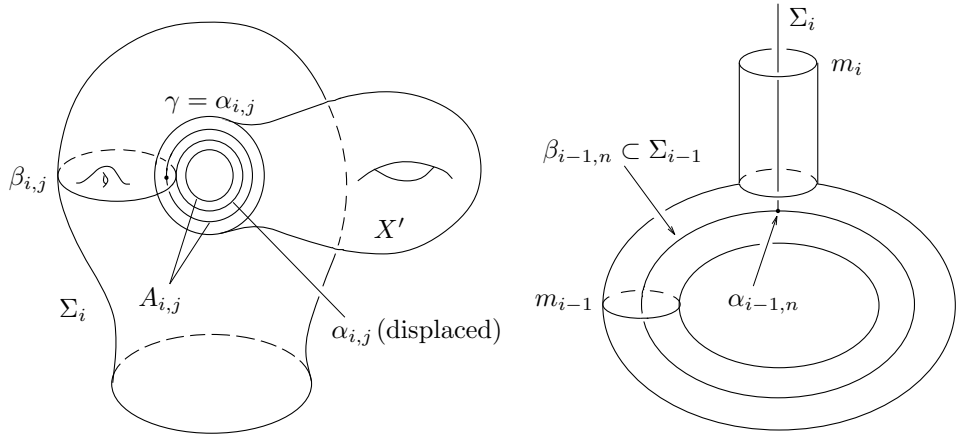
 
\torus
\caption{Steps 1 and 2} \label{torusfigure}
\end{figure}

Let $A_{i,j}$ denote the ${\epsilon}$--circle bundle
of $X_i$ in $S^4$, restricted to a parallel displacement of $\alpha_{i,j}$
in $X_i$, see Figure \ref{torusfigure}. The corresponding ${\epsilon}$--disk bundle, for $\epsilon$ small enough, can be used to see that the $2$--torus $A_{i,j}$ has linking number $1$ with $\beta_{i,j}$ 
and does not link
other curves in the collection $\{\alpha_{s,t},\beta_{s,t}\}$.  
Note that if there is a higher stage
attached to $\beta_{i,j}$ then it intersects $A_{i,j}$ in a single point,
while if there is no stage attached to $\beta_{i,j}$ then $A_{i,j}$ 
is disjoint from $X$, and the generator of 
$H_2(S^4\smallsetminus X)$ represented by $A_{i,j}$
is Alexander-dual to $\beta_{i,j}$. Similarly, let $B_{i,j}$ denote a
$2$--torus representative of the class dual to $\alpha_{i,j}$.
There are two inductive steps used in the construction of the stages of
dual Gropes.

{\em Step 1}\ \  Let $\gamma$ be a curve in the collection 
$\{{\alpha}_{i,j}, {\beta}_{i,j}\}$, and let $X'$ 
denote the subgrope of $X$ which is attached to 
$\gamma$. Since $X$ is framed 
and embedded, a parallel copy of $\gamma$ in $S^4$ bounds a parallel
copy of $X'$ in the complement of $X$. 
If there is no higher stage attached to $\gamma$ then the application
of Step 1 to this curve is empty.

{\em Step 2}\ \  Let $\Sigma_i$ be a connected 
component of the $i$-th stage of $X$, 
and let $m_i$ denote a meridian of $\Sigma_i$ 
in $S^4$, that is, $m_i$ is the boundary of a small normal disk to $\Sigma_i$
at an interior point. 
Suppose $i>1$ and let $\Sigma_{i-1}$ denote the previous stage, so that
$\Sigma_i$ is attached to $\Sigma_{i-1}$ along some curve, say 
$\alpha_{i-1,n}$. The torus $B_{i-1,n}$ meets 
$\Sigma_i$ in a point, but 
making a puncture into $B_{i-1,n}$ around this intersection point and connecting it by a tube with $m_i$
exhibits $m_i$ as the boundary of a punctured torus in the complement of $X$,
see Figure~\ref{torusfigure}. 

By construction, $H_1(X)$ is generated by those curves $\{\alpha_{i,j},\beta_{i,j}\}$
which do not have a higher stage attached to them. Fix one of these curves, 
say $\beta_{i,j}$. We will show that its dual torus $A_{i,j}$ is the first
stage of an embedded Grope $G\subset S^4\smallsetminus X$ of class $k$.
The meridian $m_i$ and a parallel copy of ${\alpha}_{i,j}$
form a symplectic basis of circles for $A_{i,j}$.
Apply Step 1 to $\alpha_{i,j}$.
If $i=1$, the result of Step 1 is a grope at least of class $k$ and we are done.
If $i>1$, apply in addition Step 2 to $m_i$. The result of Step 2 is a grope with a new
genus 1 surface stage, the tips of which are the meridian $m_{i-1}$
and a parallel copy of one of the standard curves in the previous stage, 
say $\beta_{i-1,n}$.
The next Step 1 -- Step 2 cycle is applied to these tips. Altogether
there are $i$ cycles, forming the grope $G$.

The trees corresponding to dual gropes constructed above may be
read off the tree associated to $X$, as follows. Start with the tree
$T_{X}$ for $X$, and pick the tip ($1$--valent vertex), corresponding
to the curve ${\beta}_{i,j}$. The algorithm for drawing
the tree $T_{G}$ of the grope $G$, Alexander-dual to ${\beta}_{i,j}$,
reflects Steps 1 and 2 above.  
Consider the path $p$ from ${\beta}_{i,j}$ to the root of $T_{X}$,
and start at the vertex ${\alpha}_{i-1,n}$, adjacent to ${\beta}_{i,j}$.
Erase all branches in $T_{X}$, ``growing'' from ${\alpha}_{i-1,n}$,
except for the edge $[{\beta}_{i,j}\ {\alpha}_{i-1,n}]$ which has 
been previously considered, and its
``partner'' branch $[{\alpha}_{i,j}\ {\alpha}_{i-1,n}]$, and then move one edge
down along the path $p$. This step is repeated $i$ times, until the root 
of $T_{X}$ is reached. The tree $T_{G}$ is obtained by copying the part of
$T_{X}$ which is not erased, with the tip ${\beta}_{i,j}$ drawn as
the root, see figure \ref{dualtree}.

\begin{figure}[ht]
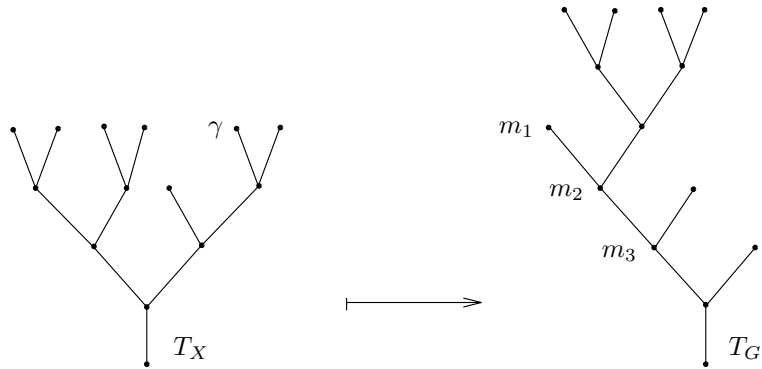
 
\genusone
\caption{Tree duality in the genus $1$ case} \label{genusone}
\end{figure}

Note the ``distinguished'' path in $T_{G}$, starting at the root and
labelled by $m_i, m_{i-1},$ $\ldots, m_1$. Each of the 
vertices $m_i, m_{i-1}, \ldots, m_2$ is trivalent (this corresponds to the
fact that all surfaces constructed by applications of Step 2 have genus $1$),
see figures \ref{dualtree}, \ref{genusone}. 
In particular, the class of $G$ may be computed
as the sum of classes of the gropes attached to the ``partner'' vertices 
of $m_i,\ldots, m_1$, plus $1$.

We will now prove that the dual grope $G$ is at least of class~$k$.
The proof is by induction on the class of $X$. For surfaces (class~$=2$)
the construction gives tori in the collection 
$\{A_{i,j}, B_{i,j}\}$.
Suppose the statement holds for Gropes of class less than $k$, and
let $X$ be a Grope of class $k$. By definition, for each standard
pair of dual circles ${\alpha}, {\beta}$ in the first stage $\Sigma$
of $X$ there is a $p$--grope $X_{\alpha}$ attached to $\alpha$ 
and a $q$--grope $X_{\beta}$ attached to $\beta$ with $p+q\geq k$.
Let $\gamma$ be one of the tips of $X_{\alpha}$. By the induction hypothesis,
the grope $G_{\alpha}$ dual to $\gamma$, given by the construction above 
for $X_{\alpha}$, is at least of class $p$. $G$ is obtained from $G_{\alpha}$
by first attaching a genus $1$ surface to $m_2$, with new tips $m_1$
and a parallel copy of $\beta$ (Step 2), and then attaching a parallel
copy of $X_{\beta}$ (Step 1). According to the computation above of the 
class of $G$ in terms of its tree, it is equal to $p+q\geq k$.

\begin{figure}[ht]
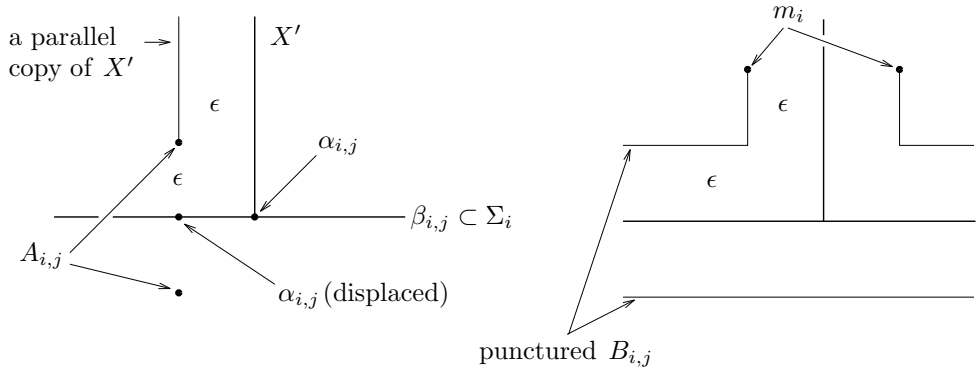
 
\epsilonfigure
\caption{Steps 1 and 2 at a distance $\epsilon$ from $X$} \label{epsilon}
\end{figure}

It remains to show that the dual gropes can be made disjoint, and that 
they are $0$--framed. Each dual grope may be arranged to lie
in the boundary of a regular $\epsilon$--neighborhood of $X$, for some
small $\epsilon$. Figure \ref{epsilon} shows how Steps 1 and 2 
are performed at a distance $\epsilon$ from $X$.
Note that although tori $A_{i,j}$ and

\eject

$B_{i,j}$ intersect, 
at most one of them is used in the construction of a dual grope for 
each index $(i,j)$. Taking distinct values 
${\epsilon}_1,\ldots, {\epsilon}_r$, the gropes are
arranged to be pairwise disjoint. The same argument shows that
each grope $G$ has a parallel copy $G'$ with $G\cap G'=\emptyset$,
hence its thickening in $S^4$ is standard.

To prove the converse part of the theorem, suppose that 
$H_2(S^4\smallsetminus X)$ is generated by $r$ disjoint maps of closed
gropes. Perturb the maps in the complement of $X$, so that they are immersions and their images have
at most a finite number of transverse self-intersection points. 
The usual pushing down and twisting procedures from \cite{FQ} produce 
closed disjoint $0$--framed gropes $G_1,\dots,G_r$ whose only failure to being
actual Gropes lies in possible self-intersections of the bottom stage surfaces.
The $G_i$ still lie in the complement of $X$ and have class~$k+1$. 
The proof of the first part of Theorem~\ref{duality} shows 
that $H_2(Y)/{\phi}_{k+1}(Y)$ is generated by the ``Clifford tori'' in the 
neighborhoods of self-intersection points of the $G_i$, 
where $Y$ denotes the complement of all $G_i$ in $S^4$.
Assume $X$ is connected (otherwise consider a connected component of $X$), 
and let $X'$ denote $X$ with a $2$--cell removed from its bottom stage.
The relations given by the Clifford tori are among the defining relations 
of the Milnor group on meridians to the gropes, and Dwyer's theorem shows
(as in \cite{FT2}, Lemma 2.6) that the inclusion map induces an isomorphism
\[ M{\pi}_1(X')/M{\pi}_1(X')^{k+1}\cong M{\pi}_1(Y)/M{\pi}_1(Y)^{k+1}. \]
\noindent
Consider the boundary curve $\gamma$ of $X'$. 
Since $X$ is a grope of class $k$, by Lemma~\ref{Milnor maximal k} below
we get ${\gamma}\notin M{\pi}_1(X')^{k+1}$. On the other hand, $\gamma$ 
bounds a disk in $Y$, hence $\gamma=1\in M{\pi}_1(Y)$. This contradiction 
concludes the proof of Theorem~\ref{duality}.
\qed

\begin{lem}\label{Milnor maximal k} \sl
Let $(G,\gamma)$ be a grope of class $k$. Then $\gamma\notin M\pi_1(G)^{k+1}$.
\end{lem}
\begin{prf}
This is best proven by an
induction on $k$, starting with the fact that
$\pi_1(\Sigma)$ is freely generated by all $\alpha_i$
and
$\beta_i$. Here $\Sigma$ is the bottom surface stage of the grope $(G,\gamma)$ with a standard symplectic basis of circles $\alpha_i, \beta_i$. The Magnus expansion for the free Milnor group (see 
\cite{Milnor1}, \cite{FT2} or the proof of Theorem~\ref{composition}) shows that $\gamma =
\prod\big[ \alpha_i, \beta_i\big]$ does not lie in
$M\pi_1(\Sigma)^3$.  Similarly, for $k>2$, $\pi_1(G)$ is
freely generated by those   circles in a standard
symplectic basis of a surface stage in $G$ to which  
nothing else is attached.  Now assume that the $k$--grope
$(G,\gamma)$ is obtained by attaching $p_i$--gropes
$G_{\alpha_i}$ to $\alpha_i$ and $q_i$--gropes
$G_{\beta_i}$ to $\beta_i,\ p_i+q_i \ge k$.  By 
induction, $\alpha_i \notin M\pi_1(G_{\alpha_i})^{
p_i+1}$ and
$\beta_i \notin M\pi_1(G_{\beta_i})^{q_i+1}$ since
$p_i, q_i
\ge 1$.  But the free generators of
$\pi_1(G_{\alpha_i})$ and $\pi_1(G_{\beta_i})$ are
contained in the set of free generators of $\pi_1(G)$
and therefore $\gamma=\prod\big[\alpha_i,
\beta_i\big] \notin M\pi_1(G)^{k+1}$.  Again, this may be
seen by applying the Magnus expansion to
$M\pi_1(G)$.
\end{prf}

\begin{rem} \rm
In the case when all stages of a Grope $X$ are tori,
the correspondence between its tree $T_{X}$ and the trees of 
the dual Gropes, given in the proof of theorem \ref{duality}, 
is particularly appealing and easy to describe. Let $\gamma$ be a tip of $T_{X}$.
The tree for the Grope, Alexander-dual to $\gamma$, is obtained
by redrawing $T_{X}$, only with $\gamma$ drawn as the root, see
Figure~\ref{genusone}.
\end{rem}

As a corollary of Theorem~\ref{duality}
we get the following result.

\begin{cor}
\label{grope concordance} \sl
Let $L=(l_1,\ldots,l_n)$ and $L'=(l'_1,\ldots,l'_n)$ be two links in
$S^3\times\{0\}$ and $S^3\times\{1\}$ respectively.
Suppose there are disjointly embedded annulus-like Gropes
$A_1,\ldots,A_n$  of class $k$ in $S^3\times [0,1]$ with
$\partial A_i=l_i\cup l'_i$, $i=1,\ldots,n$. Then there is an isomorphism of nilpotent quotients
\[ \pi_1(S^3\smallsetminus L)/\pi_1(S^3\smallsetminus L)^k\cong \pi_1(S^3\smallsetminus L')/\pi_1(S^3\smallsetminus
L')^k \]
\end{cor}

\begin{rem} \rm 
For those readers who are familiar with  Milnor's $\bar\mu$--invariants we should mention that the above statement
directly implies that
for any multi-index $I$ of length $|I|\leq k$ one gets $\bar\mu_L(I)=\bar\mu_{L'}(I)$. For a different proof of this
consequence see \cite{Krushkal}.
\end{rem}
\rk{Proof of Corollary~\ref{grope concordance}}
The proof is a $\phi_k$--version of Stallings' proof of the concordance invariance of all nilpotent quotients of
$\pi_1(S^3 \smallsetminus L)$, see \cite{Stallings}. Namely, Alexander duality and Theorem~\ref{duality} imply that the
inclusion maps

\[(S^3\times\{0\}\smallsetminus L)\hookrightarrow (S^3\times[0,1]\smallsetminus
(A_1\cup\ldots\cup A_n))\hookleftarrow (S^3\times\{1\}\smallsetminus L') \]

\noindent
induce isomorphisms on $H_1(\, .\, )$ and on $H_2(\, .\, )/\phi_k$.
So by Dwyer's Theorem they induce isomorphisms on $\pi_1/(\pi_1)^k$.
\qed

\section{The Link Composition Lemma} \label{compositionsection}

The Link Composition Lemma was originally formulated in \cite{FL}.
The reader should be cautious about its proof given there; it can
be made precise using Milnor's $\bar\mu$--invariants with repeating
coefficients, while this section presents an alternative proof.

Given a link $\widehat L =(l_1,\ldots, l_{k+1})$ in $S^3$ and
a link $Q=(q_1,\ldots, q_m)$ in the solid torus $S^1\times D^2$, 
their ``composition'' is obtained by replacing the last component
of $\widehat L$ with $Q$. More precisely, it is defined as
$L\cup{\phi}(Q)$ where $L=(l_1,\ldots, l_k)$ and 
${\phi}\co S^1\times D^2\hookrightarrow S^3$ is a $0$--framed embedding
whose image is a tubular neighborhood of $l_{k+1}$.
The meridian $\{1\}\times\partial D^2$ of the solid torus will
be denoted by $\wedge$ and we put $\widehat Q:=Q\cup\wedge$. We sometimes think of $Q$ or $\widehat Q$ as links in $S^3$ via the standard embedding $S^1 \times D^2 \hra S^3$.

\begin{figure}[ht]
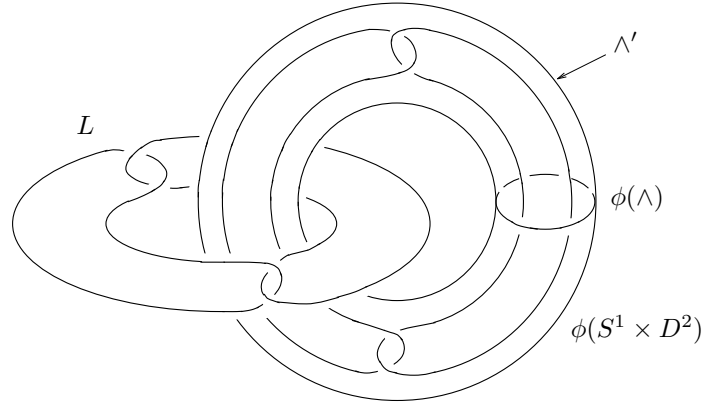
 
\composition
\caption{In this example $\widehat L$ is the Borromean rings, and 
$Q$ is the Bing double of the core circle of $S^1\times D^2$.}
\label{compositionfigure}
\end{figure}

\begin{thm}\rm (Link Composition Lemma) \label{composition} \sl \nl
{\rm(i)}\ \ If $\widehat L$ and $\widehat Q$ are both homotopically
essential in $S^3$ then $L\cup{\phi}(Q)$ is also homotopically essential.\\
{\rm(ii)}\ \ Conversely, if $L\cup{\phi}(Q)$ is homotopically essential
and if both $\widehat L$ and $\widehat Q$ are almost homotopically
trivial, then both $\widehat L$ and $\widehat Q$ are homotopically
essential in $S^3$.
\end{thm}

\begin{rem} \rm
Part (ii) does not hold without the almost triviality assumption
on $\widehat L$ and $\widehat Q$. For example, let $\widehat L$ 
consist of just one component $l_1$, and let $Q$ be a Hopf link
contained in a small ball in $S^1\times D^2$. Then $L\cup\phi(Q)={\phi}(Q)$ is
homotopically essential, yet $\widehat L$ is trivial.
\end{rem}

In part (i), if either $L$ or $Q$ is homotopically
essential, then their composition $L\cup{\phi}(Q)$ is also essential.
(Note that $\widehat Q$ and ${\phi}(\widehat Q)$ are homotopically
equivalent, see Lemma 3.2 in \cite{FL}.) 
If neither $L$ nor $Q$ is homotopically essential,
then by deleting some components of $L$ and $Q$ if necessary,
one may assume that $\widehat L$ and $\widehat Q$ are almost homotopically
trivial (and still homotopically essential).
In the case when $L\cup{\phi}(Q)$ is not almost homotopically 
trivial part (i) follows immediately. Similarly, part (ii) can be proved in this case easily by induction on the number of components of $L$ and $Q$. 

{\em Therefore, we will assume from now on that $\widehat L$, $\widehat Q$ and 
$L\cup{\phi}(Q)$ are almost homotopically trivial links in $S^3$.}

\begin{lem} \label{infgropes} \sl
If $\widehat L$ and
$\widehat Q$ are both homotopically trivial in $S^3$ then ${\phi}(\wedge)$ represents
the trivial element in the Milnor group $M(L\cup{\phi}(Q))$.
\end{lem}

\begin{prf}
Let $\wedge'$ denote ${\phi}(S^1\times\{1\})$.
The Milnor group $M(L\cup{\phi}(Q))$ is nilpotent
of class $k+m+1$, so it suffices to show that ${\phi}(\wedge)$ represents
an element in $\pi_1(S^3\smallsetminus (L\cup{\phi}(Q)))^{k+m+1}$.
This will be achieved by constructing an $\infty$--grope $G$ bounded by
${\phi}(\wedge)$ in the complement of $L\cup{\phi}(Q)$.
In fact, the construction also gives an $\infty$--grope 
$G'\subset S^3\smallsetminus (L\cup{\phi}(Q))$ bounded by $\wedge'$.

Consider $S^1\times D^2$ as a standard unknotted solid torus in $S^3$,
and let $c$ denote the core of the complementary solid torus $D^2\times S^1$.
Since $\widehat Q$ is homotopically trivial, after changing
$Q$ by an appropriate link homotopy in $S^1\times D^2$, $\wedge$ bounds

\eject 

an immersed disk $\Delta\subset S^3$ in the complement of the new link. 
Denote the new link by $Q$ again. 
Similarly $L$ can be changed so that the untwisted parallel copy $\wedge'$
of $l_{k+1}$ bounds a disk $\Delta'\subset S^3\smallsetminus L$. Recall that
$M(L\cup{\phi}(Q))$ does not change if $L\cup{\phi}(Q)$
is modified by a link homotopy.

The intersection number of $\Delta$ with $c$ is trivial, 
since $\wedge$ and $c$ do not link.
Replace the union of disks $\Delta\cap (D^2\times S^1)$ by annuli lying
in $\partial (D^2\times S^1)$ to get 
$\Sigma\subset S^1\times D^2\smallsetminus Q$, an
immersed surface bounded by $\wedge$. Similarly the intersection number
of $\Delta'$ with the core circle of ${\phi}(S^1\times D^2)$ is trivial,
and $\wedge'$ bounds $\Sigma'\subset S^3\smallsetminus (L\cup{\phi}(Q))$.
The surfaces ${\phi}(\Sigma)$ and $\Sigma'$ are the first stages
of the gropes $G$ and $G'$ respectively.

Notice that half of the basis for $H_1({\phi}(\Sigma))$ is represented
by parallel copies of $\wedge'$. They 
bound the obvious surfaces: annuli 
connecting them with $\wedge'$ union with $\Sigma'$, which
provide the second stage for $G$. Since this construction
is symmetric, it provides all higher stages for both $G$
and $G'$.
\end{prf} 

\begin{lem} \label{welldefined} \sl 
Let $i\co S^3\smallsetminus$ neighborhood $(\widehat L\smallsetminus l_1)\longrightarrow
S^3\smallsetminus (L\cup{\phi}(Q)\smallsetminus l_1)$ denote the inclusion map,
and let $i_{\#}$ be the induced map on $\pi_1$. Then $i_{\#}$ induces a well defined map $i_{*}$ 
of Milnor groups.
\end{lem}

\begin{rem} \rm
Given two groups $G$ and $H$ normally generated by $g_i$ respectively $h_j$,
let $MG$ and $MH$ be their Milnor groups defined with respect to the
given sets of normal generators. If a homomorphism
$\phi\co G\longrightarrow H$ maps each $g_i$ to one of the $h_j$
then it induces a homomorphism $M\phi\co MG\longrightarrow MH$.
In general, $\phi\co G\longrightarrow H$ induces a homomorphism
of the Milnor groups if and only if $\phi(g_i)$ commutes with
$\phi(g_i)^{\phi(g)}$ in $MH$ for all $i$ and all $g\in G$.
\end{rem}

\rk{Proof of Lemma~\ref{welldefined}}
The Milnor groups $M(\widehat L\smallsetminus l_1)$ and 
$M(L\cup{\phi}(Q)\smallsetminus l_1)$ are generated
by meridians.
Moreover,  $i_{\#}(m_i)=m_i$ for $i=2,\ldots,k$ and
$i_{\#}(m_{k+1})={\phi}(\wedge)$ where $m_1,\ldots, m_{k+1}$ are meridians to
the components of $\widehat L$.
Hence to show that $i_{*}$ is well-defined it suffices to prove
that all the commutators 
\[ [{\phi}(\wedge),({\phi}\wedge)^{i_{\#}(g)}], \quad g\in\pi_1(S^3\smallsetminus (\widehat L\smallsetminus l_1)),
\]

\noindent
are trivial in $M(L\cup{\phi}(Q)\smallsetminus l_1))$.
Consider the following exact sequence, obtained by deleting the component $q_1$ of $Q$.

\[ ker(\psi)\longrightarrow M(L\cup{\phi}(Q)\smallsetminus l_1)
\buildrel{\psi}\over\longrightarrow
M(L\cup{\phi}(Q)\smallsetminus (l_1\cup \phi(q_1)))\ra 0 \]

\eject

An application of Lemma~\ref{infgropes} to 
$(\widehat L\smallsetminus l_1)$
and to $(\widehat Q\smallsetminus q_1)$ shows that
$\psi({\phi}(\wedge))=1$ and hence 
${\phi}(\wedge),{\phi}(\wedge)^g\in ker(\psi)$.
The observation that $ker(\psi)$ is generated by the meridians
to $\phi(q_1)$ and hence is commutative finishes the proof of
Lemma~\ref{welldefined}.
\qed

\rk{Proof of Theorem~\ref{composition}}
Let $M(F_{m_1,\ldots,m_{s+1}})$ be the Milnor group of a free group, ie the Milnor group of the trivial link on $s+1$
components with meridians $m_i$. Let $R(y_1,\ldots, y_s)$
be the quotient of the free associative ring
on generators $y_1,\ldots,y_s$ by the ideal generated by the
monomials $y_{i_1}\cdots y_{i_r}$ with one index occurring at least twice. The additive group $(R(y_1,\ldots, y_s),+)$
of this ring is free abelian on generators $y_{i_1}\cdots y_{i_r}$ where all indices are distinct. Milnor
\cite{Milnor1} showed that setting $m_{s+1}=1$ induces a short exact sequence of groups
\[
1\longrightarrow (R(y_1,\ldots, y_s),+) \buildrel{r}\over\longrightarrow 
M(F_{m_1,\ldots,m_{s+1}}) \buildrel{i}\over\longrightarrow M(F_{m_1,\ldots,m_{s}})\longrightarrow 1
\]

\noindent
where $r$ is defined on the above free generators by left-iterated commutators with $m_{s+1}$:
\[
r(y_{j_1}\cdots y_{j_k}):=
[m_{j_1},[m_{j_2},\ldots,[m_{j_k},m_{s+1}]\ldots]]
\]

\noindent
In particular, $r(0)=1$ and $r(1)=m_{s+1}$. Obviously, the above extension of groups splits by sending $m_i$ to $m_i$.
This splitting induces the following conjugation action of $M(F_{m_1,\ldots,m_s})$ on $R(y_1,\ldots, y_s)$. Let
$Y:=y_{j_1}\cdots y_{j_k}$, then

\[
m_i \cdot r(Y)\cdot  m_i^{-1}  =  [m_i,r(Y)] \cdot r(Y) = \]
\[ [m_i,[m_{j_1},[m_{j_2},\ldots,[m_{j_k},m_{s+1}]\ldots]] \cdot r(Y)  =  r((y_i+1)\cdot Y)
\]

\noindent
which implies that $m_i$ acts on $R(y_1,\ldots$, $y_s)$ by ring multiplication with $y_i+1$ on the left. Since $m_i$
generate the group $M(F_{m_1,\ldots,m_s})$ this defines a well defined homomorphism of $M(F_{m_1,\ldots,m_s})$ into the
units of the ring $R(y_1,\ldots$, $y_s)$. In fact, this is the {\em Magnus expansion}, well known in the context of
free groups (rather than free Milnor groups). We conclude in particular, that the abelian group $(R(y_1,\ldots,
y_s),+)$ is generated by $y_i$ as a module over the group $M(F_{m_1,\ldots,m_s})$.

Returning to the notation of Theorem~\ref{composition}, we have the following commutative diagram of group extensions. We use the fact that the links $L\cup {\phi}(Q) \smallsetminus l_1$ and $\widehat L \smallsetminus l_1$ are
homotopically trivial. Here $y_i$ are the variables corresponding to the link $L$ and $z_j$ are the variables
corresponding to ${\phi}(Q)$. We introduce short notations
$R(\mathcal{Y}):=R(y_1,\ldots,y_k)$ and 
$R(\mathcal{Y},\mathcal{Z}):=R(y_1,\ldots,y_k,z_2,\ldots,z_m)$.

\[
\begin{CD}
R(\mathcal{Y},\mathcal{Z}) @>{r}>> M(L\cup {\phi}(Q)\smallsetminus l_1)
@>{i}>> M(L\cup {\phi}(Q)\smallsetminus (l_1\cup {\phi}(q_1))) \\
 @AA{\sigma}A @AA{lc}A @AA{lc}A   \\
R(\mathcal{Y}) @>{\bar r}>> M(\widehat L\smallsetminus l_1)
@>{j}>> M(L\smallsetminus l_1)  \\
\end{CD}
\]

\noindent
Recall that by definition $lc(m_i)=m_i$ for all meridians $m_2,\ldots,m_k$ of $L \smallsetminus l_1$. Moreover, the
link composition map $lc$ sends the meridian $m_{k+1}$ to the $\wedge$--curve of ${\phi}(Q)$.

The existence of the homomorphism $lc$ on the Milnor group level already implies our claim~(ii) in Theorem~\ref{composition}: By assumption,
$l_1$ represents the trivial element in $M(\widehat L\smallsetminus l_1)$ since $\widehat L$ is homotopically trivial.
Consequently, $lc(l_1)=l_1$ is also trivial in $M(L\cup {\phi}(Q)\smallsetminus l_1)$ and hence by \cite{Milnor1} the link
$L\cup {\phi}(Q)$ is homotopically trivial.

The key fact in our approach to part~(i) of Theorem~\ref{composition} is the following result which says that link
composition corresponds to ring multiplication.
\begin{lem} \label{multiplication} \sl
The homomorphism $\sigma\co \! R(y_2,\ldots,y_{k})\longrightarrow
R(y_2,\ldots,y_{k}, z_2,\ldots,z_m)$ is given by ring multiplication
with $r^{-1}(\wedge)$ on the right.
\end{lem}
Note that by Lemma~\ref{infgropes} $\wedge$ is trivial in
$M(L\cup {\phi}(Q)\smallsetminus (l_1\cup {\phi}(q_1)))$, so that it makes sense
to consider $r^{-1}(\wedge)$. We will abbreviate this important element by $\wedge_R$.

\rk{Proof of Lemma~\ref{multiplication}}
Since the above diagram commutes and $(R(y_2,\ldots, y_{k}),$ $+)$ is generated by $y_i$ as a module over the group
$M(F_{m_2,\ldots,m_{k}})$ it suffices to check our claim for these generators $y_i$. We get by definition

\[
lc(\bar r(y_i)) = lc([m_i,m_{k+1}]) =[m_i,\wedge]= (m_i\cdot\wedge\cdot m_i^{-1})\cdot\wedge^{-1} \]
\[ = r((y_i+1)\cdot\wedge_R)\cdot \wedge^{-1} 
 =  r(y_i\cdot\wedge_R).
\]

\noindent
We are using the fact that conjugation by $m_i$ corresponds to left multiplication by $(y_i+1)$.
\qed

Since $L$ is homotopically trivial and $\widehat L$ is homotopically
essential, it follows that  $0\neq l_1\in ker(j)$. After possibly reordering the $y_i$ this implies in addition that
for some integer $a\neq 0$ we have

\[
\bar r^{-1}(l_1)=a\cdot (y_2\cdots y_k) +
{\rm \; terms \; obtained \; by \; permutations \; from \; } y_2\cdots y_k.
\]

\noindent
Setting all the meridians $m_i$ of $L$ to $1$ (which implies setting the variables $y_i$ to $0$), we get a commutative
diagram of group extensions

\[
\begin{CD}
 R(\mathcal{Z}) @>r>> M({\phi}(Q))
@>{i}>> M({\phi}(Q)\smallsetminus {\phi}(q_1)) \\
@AA{p}A @AAA @AAA   \\
R(\mathcal{Y},\mathcal{Z}) @>r>> M(L\cup {\phi}(Q)\smallsetminus l_1)
@>{i}>> M(L\cup {\phi}(Q)\smallsetminus (l_1\cup {\phi}(q_1)))\\
\end{CD}
\]

\noindent
As before, $R(\mathcal{Z})$ and $R(\mathcal{Y},\mathcal{Z})$
are short notations for $R(z_2,\ldots,z_m)$ and 
$R(y_1,$ $\ldots,y_k,z_2,\ldots,z_m)$ respectively.
Since $\widehat Q$  (and equivalently $\phi(\widehat Q)$) is homotopically essential
we have $0\neq \wedge\in ker(i)$.
This shows that $p(\wedge_R)\neq 0$. The almost triviality of $\widehat Q$ implies in addition that after possibly
reordering the $z_j$ we have for some integer $b\neq 0$

\[
p(\wedge_R)=b\cdot (z_2\cdots z_m) + {\rm \; terms \; obtained \; 
by \; permutations \; from \; } z_2\cdots z_m.
\]

\noindent
It follows from Lemma~\ref{multiplication} that
$r^{-1}(l_1)=\bar r^{-1}(l_1)\cdot
\wedge_R$. This product contains the term
\[ a b \cdot
(y_2\cdots y_{k}\cdot z_2\cdots z_m),
 \]

\noindent
the coefficient $ab$ of which is non-zero. This completes the proof of
Theorem~\ref{composition}.\break\hbox{}
\qed

\begin{rem} \rm
Those readers who are familiar with Milnor's $\bar\mu$--invariants will have recognized that the above proof in fact
shows that the first non-vanishing $\bar\mu$--invariants are multiplicative under link composition.
\end{rem}

\vfil\eject

\end{document}